
\def\DJ{\leavevmode\setbox0=\hbox{D}\kern0pt\rlap
{\kern.04em\raise.188\ht0\hbox{-}}D}
\nopagenumbers
\font\aa=cmss12 at 12pt
\font\bb=cmcsc10

\font\ff=cmr9
\font\jj=cmr8
\font\mm=cmti8

\font\rr=cmcsc10 at 12pt

\footnote{{}}{\hskip -2 mm  {\jj 2001}     {\mm   Mathematics Subject
Classification}. {\jj Primary 11F72, Secondary 11M06}.

\hskip -2 mm
{\mm Key words and phrases: } {\jj Fourier coefficients of cusp forms,
Ramanujan-Petersson conjecture, Riemann zeta-function}}


\baselineskip=13pt
\def\hf{{\textstyle{1\over2}}}
\def\a{\alpha}\def\b{\beta}

\def\e{\varepsilon}
\def\f{\varphi}
\def\G{\Gamma}
\def\k{\kappa_j}
\def\s{\sigma}

\def\={\;=\;}

\def\no{\noindent}
\def\R{\Re{\rm e}\,}  \def\s{\sigma}
\def\z{\zeta}

\def\no{\noindent} 
  \def\={\,=\,}
\def\hf{{\textstyle{1\over2}}}

\def\f{\varphi}

\font\tenmsb=msbm10
\font\sevenmsb=msbm7
\font\fivemsb=msbm5
\newfam\msbfam
\textfont\msbfam=\tenmsb
\scriptfont\msbfam=\sevenmsb
\scriptscriptfont\msbfam=\fivemsb
\def\Bbb#1{{\fam\msbfam #1}}

\def \NN {\Bbb N}

\def \RR {\Bbb R}
\def \ZZ {\Bbb Z}

\def\rightheadline{\hfil\ff  On sums of Fourier coefficients
of cusp forms\hfil\folio}
\def\leftheadline{\ff\folio\hfil Aleksandar Ivi\'c
 \hfil}
\def\emptyheadline{\hfil}
\headline{\ifnum\pageno=1 \emptyheadline\else
\ifodd\pageno \rightheadline \else \leftheadline\fi\fi}
\topglue2cm
\centerline{\aa  ON SUMS  OF  FOURIER COEFFICIENTS OF CUSP FORMS }

\bigskip\bigskip
\centerline{\rr Aleksandar Ivi\'c}
\bigskip\bigskip

The aim of this note is to evaluate asymptotically the sum
$$
F(x) \;:=\; \sum_{n\le x} f(n^2)\leqno(1)
$$
in case $f(n)$ is the Fourier coefficient of a holomorphic or
non-holomorphic cusp form. We shall first deal with the latter
case, which is more complicated. Let  as usual
$\,\{\lambda_j = \k^2 + {1\over4}\} \,\cup\, \{0\}\,$ be the
discrete spectrum of the non-Euclidean Laplacian acting on
$\,SL(2,\ZZ)\,$--automorphic forms. Further let $\rho_j(n)$ denote
the $n$-th Fourier coefficient of the Maass wave form $\f_j(z)$
corresponding to the eigenvalue $\lambda_j$ to which the Hecke series
$$
H_j(s) \= \sum_{n=1}^\infty t_j(n)n^{-s}\qquad(\R s > 1)
$$
is attached (see Y. Motohashi [6] for an extensive account). For
every $n \in \NN$ we have
$\rho_j(n) = \rho_j(1)t_j(n)$, so that one may consider sums of $t_j(n)$
instead of sums of $\rho_j(n)$. In [4] T. Meurman and the author proved
that, for $\k \le x^{1-\a}$, we have uniformly in $\k$
$$
\sum_{n\le x}t_j^2(n) \= {12x\over\pi^2\a_j} + O_\e(x^\e R(x))\leqno(2)
$$
with
$$
R(x) \;:=\; \k^{1\over2-2\a}x^{1\over2} +
\min\left(\k^{1+10\a\over3+6\a}x^{1\over2} + x^{3+6\a\over5},\,
\k^{-1}x^{1+2\a}\right).\leqno(3)
$$
Here as usual
$$
\a_j \= {|\rho_j(1)|^2\over\cosh(\pi\k)},\leqno(4)
$$
and $\a \;(\ge 0)$ is the constant for which
$$
t_j(n) \;\ll_\e\; n^{\a+\e}\leqno(5)
$$
holds uniformly in $\k$. Moreover, $\e$ denotes arbitrarily small
positive constants, not necessarily the same ones at  each occurrence,
while $\ll_\e$ means that the $\ll$-constant depends on $\e$.
The accent in (2)-(3) is on uniformity in $\k$, since in many
applications $\k$ may vary with $x$. It is known that
(5) holds with $\a \le {5\over28}$
(see D. Bump et al. [1]). In what concerns the
order of $\a_j$, we have
$$
\k^{-\e} \;\ll_\e\; \a_j \;\ll_\e\;\k^{\e}.\leqno(6)
$$
The lower bound in (6) was proved by H. Iwaniec [5], and the upper
bound by Hoffstein--Lockhart [2].

\medskip
It seems that there are no upper bounds in the literature for the sum
$F(x)$ in (1) when $f(n) = t_j(n)$. We shall prove

\medskip
THEOREM 1. {\it For $\k\le x^c$,}
$$
0 \;<\; c \;<\; \min\left({3+6\a\over2+20\a},\,1-\a\right)\leqno(7)
$$
{\it and a suitable constant $A>0$ we have, uniformly in $\k$,}
$$
\sum_{n\le x}t_j(n^2) \;\ll\;\a_j^{-1}x
\exp\left(-A\log^{3/5}x(\log\log x)^{-1/5}\right).\leqno(8)
$$

\bigskip
{\bf Corollary 1}. With the value $\a \le {5\over28}$ it follows that
(8) holds uniformly in $\k$ for $\k\le x^c$, and any constant $c$
satisfying $0 < c < {57\over78}$.
\bigskip
{\bf Corollary 2}. If the Ramanujan-Petersson conjecture
that $\a = 0$ is true, then (8)
holds uniformly in $\k$ for $\k \le x^c$, and any constant $0<c<1$.

\bigskip
{\bf Remark 1}. The negative exponential in (8) comes from the
sharpest known error term in the prime number theorem in the form
(see e.g., [3, Chapter 12])
$$
\sum_{n\le x}\mu(n) \;\ll\;
x\exp\left(-C\log^{3/5}x(\log\log x)^{-1/5}\right)\qquad(C>0).\leqno(9)
$$
Sharper forms of the prime number theorem, which would follow from a
better zero-free region for the Riemann zeta-function $\zeta(s)$
than the one that is currently known
(see [3, Chapter 6]), would therefore lead to a better estimate than (8).

\bigskip
{\bf Remark 2}. By a result of T. Meurman and the author [4] one has
$$
\sum_{n\le x} t_j(n) \;\ll_\e\; \k^{1+\e}
$$
uniformly for $\sqrt{x} < \k \le x$, which may be compared to the bound
in (8).

\bigskip
{\bf Remark 3}. The oscillatory nature of the function $t_j(n)$ accounts
for the lack of a main term in (8). However, if one looks at the
problem of evaluating $F(x)$ in (1) when $f(n) = d(n)$,
the number of divisors of $n$, then there will be
a main term in the corresponding formula for the summatory function.
Namely the function $d(n^2)$ is generated by $\z^3(s)/\z(2s)$, which
has a pole of order three at $s=1$.
Consequently we have (see [3, eq. (14.29)]), for  suitable constants
$B_1 \,(> 0), B_2, B_3,  C\,( >0)$,
$$
\sum_{n\le x}d(n^2) = x(B_1\log^2x + B_2\log x + B_3)
+ O\left(x\exp\left(-C\log^{3/5}x(\log\log x)^{-1/5}\right)\right).
$$
\bigskip
{\bf Proof of Theorem 1}. From the multiplicative property
(see [6, eq. (3.2.8)])
$$
t_j(mn) \= \sum_{d|(m,n)}\mu(d)t_j\left({m\over d}\right)
t_j\left({n\over d}\right),
$$
where $\mu(n)$ is the M\"obius function, one has
$$
t_j(n^2) \= \sum_{d|n}\mu(d)t_j^2\left({n\over d}\right).\leqno(10)
$$
Therefore (10) gives
$$
\eqalign{
&\sum_{n\le x}t_j(n^2) = \sum_{mn\le x}\mu(m)t_j^2(n)\cr&
= \sum_{m\le\sqrt{x}}\mu(m)\sum_{n\le x/m}t_j^2(n)
+ \sum_{n\le\sqrt{x}}t_j(n^2)\sum_{\sqrt{x}<m\le x/n}\mu(m) \cr&
= \sum\nolimits_1 +\sum\nolimits_2,
\cr}
$$
say. We set for brevity
$$
\eta(x) \= (\log x)^{3/5}(\log\log x)^{-1/5},
$$
and let $0 < \b < 1,\,\b = \b(\a)$  be such a number for which (2)-(3)
give
$$
\sum_{n\le x}t_j^2(n) = {12x\over\pi^2\a_j} + O(x^\b)\leqno(11)
$$
uniformly for $\k \le x^c$. If $C$ denotes generic positive constants, then
$$
\eqalign{
\sum\nolimits_1 &= {12x\over\pi^2\a_j}\sum_{m\le\sqrt{x}}{\mu(m)\over m}
+ O\Bigl(\sum_{m\le\sqrt{x}}\left({x\over m}\right)^\b\Bigr)\cr&
= {12x\over\pi^2\a_j}\sum_{m>\sqrt{x}}{\mu(m)\over m}
+ O(x^\b x^{1-\b\over2}) \cr&
\ll \a_j^{-1}x{\rm e}^{-C\eta(x)},
\cr}
$$
since in view of $0 < \b < 1$ and the upper bound in  (6) we have
$$
x^{1+\b\over2} \;\ll\; \a_j^{-1}x{\rm e}^{-C\eta(x)}.
$$
Here we also used (9), partial summation and the well-known fact that
$$
\sum_{m=1}^\infty {\mu(m)\over m} \= 0
$$
to deduce that
$$
\sum_{m>\sqrt{x}}{\mu(m)\over m} \;\ll \;\exp(-C\eta(x)).\leqno(12)
$$
We also have, on using (9) and (11),
$$
\eqalign{
\sum\nolimits_2 &\ll \sum_{n\le\sqrt{x}}t_j^2(n)\left({x\over n}
{\rm e}^{-C\eta(x/n)} + \sqrt{x}{\rm e}^{-C\eta(x)}\right)\cr&
\ll {\rm e}^{-C\eta(x)}\Bigl(x\sum_{n\le\sqrt{x}}{t_j^2(n)\over n}
+ \sqrt{x}\sum_{n\le\sqrt{x}}t_j^2(n)\Bigr)\cr&
\ll \a_j^{-1}x{\rm e}^{-C\eta(x)}.
\cr}
$$
Now note that for $\k \le x^c,\,c < 1 - \a$ we have
$$
\k^{1\over2-2\a}x^{1\over2} \;\le\; x^\b,\quad
\b = {c\over2-2\a}+{1\over2}  \quad(\,<1).
$$
Moreover
$$
{1+10\a\over3+6\a}c \;<\;{1\over2}\quad{\rm{for}}\quad
c < {3+6\a\over2+20\a},
$$
and $\,{3+6\a\over5} < 1$. Hence (11) is satisfied if $\k \le x^c$
and $c$ is any constant satisfying (7). This finishes the proof
of Theorem 1.

\bigskip
The foregoing analysis may be applied to the case of a
holomorphic cusp form
$$
\f(z) \= \sum_{n=1}^\infty a(n){\rm e}^{2\pi inz}
$$
of weight $\kappa$ with respect  to $SL(2,\,\ZZ)$. If $\f(z)$ is
a normalized eigenform, i.e., an eigenfunction with respect to all
Hecke operators and satisfies $a(1) = 1$, then all $a(n) \in \RR$.
We have (see e.g., R.A. Rankin [7] and [8])
$$
\sum_{n=1}^\infty a(n)n^{-s} = \prod_p(1-a(p)p^{-s}+p^{\kappa-1-2s})^{-1}
\quad(\s > \hf(\kappa + 1)),
$$
and
$$
a(m)a(n) \= \sum_{d|(m,n)}d^{\kappa-1}a\left({mn\over d^2}\right)
\qquad(m,n \in \NN).\leqno(15)
$$
If we introduce the ``normalized" function $\tilde a(n)$, namely
$$
a(n) \= \tilde a(n)n^{\kappa-1\over2},\leqno(16)
$$
then we can write (15) as
$$
\tilde a(m)\tilde a(n) \= \sum_{d|(m,n)}\tilde a\left({mn\over d^2}\right)
\qquad(m,n \in \NN).\leqno(17)
$$
When $m = n$, (17) gives by the M\"obius inversion formula
$$
\tilde a(n^2) \= \sum_{d|n}\mu(d)\left(\tilde a\bigl({n\over d}
\bigr)\right)^2,\leqno(18)
$$
which is analogous to (10). If we use (18), the Rankin-Selberg formula
$$
\sum_{n\le x}\tilde a^2(n) \= Cx + O(x^{3/5})\qquad(C > 0),
$$
then by employing the method of proof of Theorem 1 we shall obtain

\bigskip
THEOREM 2. {\it We have}
$$
\sum_{n\le x}\tilde a(n^2) \;\ll\; x\exp\left(-C(\log x)^{3/5}
(\log\log x)^{-1/5}\right)\qquad(C>0).
$$

\bigskip\no
From Theorem 2 we obtain by partial summation, on using (16),

\bigskip
{\bf Corollary 3}. {\it We have}
$$
\sum_{n\le x} a(n^2) \;\ll\; x^\kappa\exp\left(-C(\log x)^{3/5}
(\log\log x)^{-1/5}\right)\qquad(C>0).
$$

\vfill\break
\topglue2cm
\centerline{\rr References}
\bigskip
\item{[1]} D. Bump, W. Duke, J. Hoffstein and H. Iwaniec, An estimate
for Fourier coefficients of Maass wave forms, Duke Math. J.
{\bf66}(1992), 75-81.

\item{[2]} J. Hoffstein and P. Lockhart, Coefficients of Maass
forms and the Siegel zero, Annals of Math. {\bf140}(1994), 161-176.

\item{[3]} A. Ivi\'c, The Riemann zeta-function, John Wiley \& Sons,
New York, 1985.

\item{[4]} A. Ivi\'c and T. Meurman, Sums of coefficients of
Hecke series, Acta Arith. {\bf68}(1994), 341-368.

\item{[5]} H. Iwaniec, Small eigenvalues of the Laplacian for
$\G_0(N)$, Acta Arith. \ {\bf56}\ (1990),\ 65-82.

\item {[6]} Y. Motohashi,  Spectral theory of the Riemann
zeta-function,  Cambridge University Press, Cambridge, 1997.

\item{[7]} R.A. Rankin, Modular forms and functions,
Cambridge University Press, Cambridge, 1977.

\item{[8]} R.A. Rankin, Fourier coefficients of cusp forms,
Math. Proc. Camb. Phil. Soc. {\bf100}(1986), 5-29.

\vskip3cm
\leftline{\bb Aleksandar Ivi\'c}
\leftline{\bb Katedra Matematike RGF-a}
\leftline{\bb Universiteta u Beogradu, \DJ u\v sina 7,}
\leftline{\bb 11000 Beograd, Serbia (Yugoslavia)}
\leftline{e-mail: {\tt aivic@matf.bg.ac.yu, aivic@rgf.bg.ac.yu}}
\bye